# The number of countable models of a countable supersimple theory

Byunghan Kim

October 1, 1996

**Abstract**

In this paper, we prove the number of countable models of a countable supersimple theory is either 1 or infinite. This result is an extension of Lachlan's theorem on a superstable theory.

## 1 Introduction

The aim of this paper is to extend the following classical result of Lachlan in a supersimple theory context.

**Theorem 1.1** *Let $T$ be a countable superstable theory. Then the number of (nonisomorphic) countable models of $T$ is either $1$ or $\geq \aleph_0$.*

In other words, we will prove the following.

**Theorem 1.2** *Let $T$ be a countable supersimple theory. Then the number of countable models of $T$ is either $1$ or $\geq \aleph_0$.*

Let us quote a brief history of Lachlan's Theorem from Baldwin ([1, XIII.2.31]).

> 'Theorem 1.1' has a long history. The first step in this direction is the proof by Baldwin and Lachlan ([2]) that the conclusion holds for countable theories which are $\aleph_1$-categorical. This proof used many special properties of $\aleph_1$-categorical theories... Then Lachlan ([6]) proved 'Theorem 1.1' by a complicated argument using rank. Lascar ([7]) simplified the proof by the use of $U$-rank... Finally, Pillay ([9]) has given an even simpler proof...

Lascar's proof of Lachlan's Theorem is essentially using the characteristics of "weight". Pillay's proof, according to a personal conversation with him, is actually a translation of Lachlan's original proof into forking context. Pillay's proof only uses the basic properties of forking ( for example, the notion of weight is not used), together with the Open Map Theorem. However as the Open Map Theorem is no longer true in a simple unstable theory, we are not able to copy the same proof for Theorem 1.2.

**Example 1** Let $M$ be the countable bipartite random graph, consisting of disjoint infinite sets $U, V$ with the relation $R$ between $U, V$. Hence for any finite disjoint subsets $X, Y$ of $U$, there is $z \in V$ such that $xRz$ for $x \in X$ and $\neg yRz$ for $y \in Y$, and vice versa. Let $A = \{a_i | i < \omega\} \subseteq U$. Choose $c \in U \setminus A$ so that $tp(c/A)$ is not isolated. Also select $b \in V$ such that $\neg a_i R b$ for all $i$, and $cRb$. Then $tp(c/Ab)$ does not fork over $A$, whereas $tp(c/Ab)$ is isolated.

In fact, one can come up with the following version of the Open Map Theorem for a simple theory, using the exactly same proof of the Open Map Theorem for a stable theory with Fact 1.4 (see the proof of 4.27 in [8]). But the following theorem will not be used in this paper.

**Theorem 1.3** *Let $T$ be simple, and let $A \subseteq B$. For each formula $\varphi(\bar{x}) \in L(B)$, there is a (partial) type $\Delta(\bar{x})$ over $A$ such that, for each $p \in S(A)$, $\Delta(\bar{x}) \subseteq p$ iff $\varphi(\bar{x})$ is in some nonforking extension of $p$.*

**Fact 1.4** *([10]) Let $T$ be simple and let $p(\bar{x}) \in S(A)$. For each L-formula $\varphi(\bar{x}, \bar{y})$, there is a corresponding partial type $\Delta(\bar{y})$ over $A$ such that, for any $\bar{c}, \models \Delta(\bar{c})$ iff $\varphi(\bar{x}, \bar{c})$ is in some nonforking extension of $p$.*

The main novelty of our argument here is that we find a new proof of Lachlan's Theorem which uses only the symmetry and transitivity of nonforking. Hence this proof also works for Theorem 1.2.

Now we recall from [4], [5], [11], some basic facts and definitions we need. A type $p$ *forks* over a set $A$, if there are an $L$-formula $\varphi(\bar{x}, \bar{y})$ and a set of tuples $\{\bar{c}_i | i < \omega\}$ such that $p \vdash \varphi(\bar{x}, \bar{c}_0)$, $tp(\bar{c}_i/A) = tp(\bar{c}_0/A)$ for all $i < \omega$, and $\{\varphi(\bar{x}, \bar{c}_i) | i < \omega\}$ is $k$-inconsistent for some $k \in \omega$. A first order complete theory $T$ is said to be *simple* if, for any type $p \in S(B)$, $p$ does not fork over some subset $A$ of $B$ with $|A| \leq |T|$. The theory $T$ is called *supersimple* if, for any type $p \in S(B)$, $p$ does not fork over some finite subset $A$ of $B$. Hence obviously a supersimple theory is simple. Moreover $T$ is supersimple if and only if there do not exist $A_0 \subseteq A_1 \subseteq ... \subseteq A_i...$ and $p_i \in S(A_i)$ for $i < \omega$, such that $p_{i+1}$ is a forking extension of $p_i$ for each $i < \omega$. We also recall that $T$ is *unstable* if there are a formula $\psi(\bar{x}, \bar{y})$ and tuples $\bar{b}_i, \bar{c}_i$ ($i < \omega$) such that $\models \psi(\bar{b}_i, \bar{c}_j)$ iff $i \leq j \in \omega$. A theory $T$ is said to be *stable* if $T$ is not unstable, and *superstable* if $T$ is stable and supersimple. Every stable theory is simple.

In [4], it is shown that, for simple $T$, nonforking satisfies (i) *extension* : for any $p \in S(A)$ and $A \subseteq B$, $p$ has a nonforking $q$ in $S(B)$, (ii) *symmetry*: $tp(\bar{b}/A\bar{c})$ does not fork over $A$ iff $tp(\bar{c}/A\bar{b})$ does not fork over $A$, and (iii) *transitivity*: if $A \subseteq B \subseteq C$ and $p \in S(C)$, then $p$ does not fork over $A$ iff $p$ does not fork over $B$ and the restriction of $p$ to $B$ does not fork over $A$. Hence nonforking supplies a nice notion of independence to an arbitrary simple theory. If $T$ is simple, we say $\{C_i | i \in I\}$ is *independent over* $A$ if for each $i \in I$ and $\bar{c} \in C_i$, $tp(\bar{c}/A \cup \bigcup\{C_j | j \neq i, j \in I\})$ does not fork over $A$.

One of the important properties of nonforking in a simple theory is the so called *Independence Theorem*, which is not so relevant to this paper, but worth while to mention. Tuples $\bar{a}, \bar{b}$ are said to have the same *Lascar strong type over* $A$ ($Lstp(\bar{a}/A) = Lstp(\bar{b}/A)$) if there are models $M_1, ..., M_k$, each of which contains $A$, and tuples $\bar{a} = \bar{a}_0, \bar{a}_1, ..., \bar{a}_k = \bar{b}$ such that $tp(\bar{a}_{i-1}/M_i) = tp(\bar{a}_i/M_i)$ for $1 \leq i \leq k$. In [5], the following is shown.

**Fact 1.5** *(The Independence Theorem for Lascar strong types) Assume that $T$ is simple. Let $\{B, C\}$ is independent (over $\phi$). If $Lstp(\bar{d}) = Lstp(\bar{e})$ and $tp(\bar{d}/B), tp(\bar{e}/C)$ both do not fork over $\phi$, then there is $\bar{a}$ such that $tp(\bar{e}/C) \cup tp(\bar{d}/B) \subseteq tp(\bar{a}/BC)$, $tp(\bar{a}/BC)$ does not fork over $\phi$, and $Lstp(\bar{a}) = Lstp(\bar{d})$.*

The notation here is fairly standard. $T$ is a complete theory with no finite models in a first order language $L$. Types, denoted by $p, q$, are $n$-types and possibly partial. We fix a huge $\bar{\kappa}$-saturated model $\bar{M}$, as usual. Tuples $\bar{a}, \bar{b}, \bar{c}... \in \bar{M}$ are finite. Sets $A, B, C...$ are subsets of $\bar{M}$ and models which we mention are elementary submodels of $\bar{M}$, the cardinalities of all of those are strictly less than $\bar{\kappa}$.

## 2 Forking and isolation

Let us recall Pillay's notion of semi-isolation ([3, §2],[9]). We say $tp(\bar{b}/\bar{a})$ is semi-isolated if there is a formula $\varphi(\bar{x}, \bar{a})$ in $tp(\bar{b}/\bar{a})$ such that $\models \varphi(\bar{x}, \bar{a}) \to tp(\bar{b})$. Definition implies the following easy, but important facts.

**Fact 2.1** *(i) If $tp(\bar{b}/\bar{a})$ is isolated, then $tp(\bar{b}/\bar{a})$ is semi-isolated.*
  *(ii) If $tp(\bar{c}/\bar{b})$ and $tp(\bar{b}/\bar{a})$ are semi-isolated, then $tp(\bar{c}/\bar{a})$ is semi-isolated.*

**Example 2** (i) The notions semi-isolation and isolation are different. For consider the model $(Z, S)$, where $S$ is the successor function. If $a, b$ are in different chains, then $tp(b/a)$ is not isolated, but semi-isolated.
  (ii) Let $L = \{E_i | i < \omega\}$. Let $T$ be a theory saying that $E_0$ is an equivalence relation having two infinite classes, and for each $i < \omega$, equivalence relation $E_{i+1}$ refines every $E_i$-class into exactly two infinite $E_{i+1}$-classes. Then $T$ is superstable. Now if $a, c$ are in the same $E_i$-class for each $i$, and $\neg bE_0 a$, then $tp(c/b), tp(b/a)$ are isolated, while $tp(c/a)$ is not isolated ( but semi-isolated).

**Fact 2.2** *Suppose that $tp(\bar{b}/\bar{a})$ is isolated, whereas $tp(\bar{a}/\bar{b})$ is nonisolated. Then $tp(\bar{a}/\bar{b})$ is nonsemi-isolated.*

*Proof.* Suppose that $\varphi(\bar{x}, \bar{a})$ isolates $tp(\bar{b}/\bar{a})$. In order to induce a contradiction, assume that $tp(\bar{a}/\bar{b})$ is semi-isolated witnessed by $\psi(\bar{b}, \bar{y})$. Now as $tp(\bar{a}/\bar{b})$ is nonisolated, there is a formula $\phi(\bar{x}, \bar{y}) \in L$ such that $\varphi(\bar{b}, \bar{y}) \wedge \psi(\bar{b}, \bar{y}) \wedge \phi(\bar{b}, \bar{y})$ and $\varphi(\bar{b}, \bar{y}) \wedge \psi(\bar{b}, \bar{y}) \wedge \neg \phi(\bar{b}, \bar{y})$ are both consistent. Moreover both formulas imply $tp(\bar{a})$. Hence $\varphi(\bar{x}, \bar{a}) \wedge \phi(\bar{x}, \bar{a})$ and $\varphi(\bar{x}, \bar{a}) \wedge \neg \phi(\bar{x}, \bar{a})$ are both consistent. This contradicts the fact that $\varphi(\bar{x}, \bar{a})$ is a principal formula. □

Now we state a key proposition which describes the relationship between isolation and forking in a simple theory.

**Proposition 2.3** *Assume that $T$ is simple. Let $\bar{a}, \bar{b}$ be two realizations of a complete type over $\phi$. If $tp(\bar{b}/\bar{a})$ is semi-isolated, and $tp(\bar{a}/\bar{b})$ is nonsemi-isolated, then $tp(\bar{a}/\bar{b})$ forks over $\phi$.*

*Proof.* Suppose that $tp(\bar{b}/\bar{a})$ is semi-isolated witnessed by $\varphi(\bar{x}, \bar{a})$. Let $\bar{c}$ be any tuple such that $tp(\bar{c}\bar{b}) = tp(\bar{b}\bar{a})$. We claim that $\varphi(\bar{c}, \bar{x}) \wedge \varphi(\bar{x}, \bar{a})$ forks over $\phi$, in other words, $tp(\bar{b}/\bar{c}\bar{a})$ forks over $\phi$: First, let $\bar{c}_0 = \bar{c}, \bar{b}_0 = \bar{b}, \bar{a}_0 = \bar{a}$. Now there is a sequence of tuples $\langle \bar{c}_i \bar{b}_i \bar{a}_i | i < \omega \rangle$ such that, for all $i < \omega$, $tp(\bar{c}_i \bar{b}_i \bar{a}_i) = tp(\bar{c}\bar{b}\bar{a})$ and $tp(\bar{a}_{i+1} \bar{c}_i) = tp(\bar{b}\bar{a})$. We note that, by Fact 2.1, $tp(\bar{a}_j/\bar{a}_i)$ is semi-isolated for every $j \geq i$ (*). It suffices to show $\{\varphi(\bar{c}_i, \bar{x}) \wedge \varphi(\bar{x}, \bar{a}_i) | i < \omega\}$ is 2-inconsistent. If it were not 2-inconsistent, then there is $\bar{d}$ such that $\varphi(\bar{d}, \bar{a}_j)$ and $\varphi(\bar{c}_i, \bar{d})$ for some $j > i$. Therefore clearly $tp(\bar{d}/\bar{a}_j), tp(\bar{c}_i/\bar{d})$ are both semi-isolated, and hence again by Fact 2.1, so does $tp(\bar{c}_i/\bar{a}_j)$. Now as $tp(\bar{a}_j/\bar{a}_{i+1})$ is semi-isolated (by (*)), once more Fact 2.1 implies $tp(\bar{c}_i/\bar{a}_{i+1})$ is semi-isolated. But since $tp(\bar{c}_i \bar{a}_{i+1}) = tp(\bar{a}\bar{b})$, it leads a contradiction. Hence the claim is proved.

Now if $\{\bar{a}, \bar{b}\}$ is independent (over $\phi$), then by the extension, symmetry and transitivity of nonforking, we can find a tuple $\bar{c}'$ such that $tp(\bar{c}'\bar{b}) = tp(\bar{b}\bar{a})$ and $\{\bar{a}, \bar{b}, \bar{c}'\}$ is independent. This contradicts the claim above. Thus $tp(\bar{a}/\bar{b})$ forks over $\phi$. □

**Corollary 2.4** *Assume that $T$ is simple. Let $\bar{a}, \bar{b}$ be realizations of a complete type over $\phi$. If $tp(\bar{b}/\bar{a})$ is isolated, and $tp(\bar{a}/\bar{b})$ is nonisolated, then $tp(\bar{a}/\bar{b})$ forks over $\phi$.*

**Remark 2.5** (i) The simplicity of $T$ is essential in Proposition 2.3. Let $(M, <, \{c_i\}_{i<\omega})$ be the Ehrenfeucht model having 3 nonisomorphic models. The theory of the model is not simple. Choose $a, b$ such that $c_i < a < b$ for all $i$. Then $tp(a) = tp(b)$, and $tp(b/a)$ is isolated whereas $tp(a/b)$ is not isolated. But $tp(a/b), tp(b/a)$ both do not fork over $\phi$. In fact, whenever $tp(ab) = tp(bc)$, then $tp(b/ac)$ forks over $\phi$.

(ii) In 2.3, the Independence Theorem for Lascar strong types yields a cheap proof, provided there is an additional assumption that $Lstp(\bar{a}) = Lstp(\bar{b})$. Now if $\{\bar{a}, \bar{b}\}$ were independent, then there is a common realization $\bar{c}$ of $tp(\bar{d}/\bar{a})$ and $tp(\bar{e}/\bar{b})$ where $tp(\bar{d}\bar{a}) = tp(\bar{a}\bar{b}) = tp(\bar{b}\bar{e})$ and $Lstp(\bar{d}) = Lstp(\bar{a}) = Lstp(\bar{e})$. We note that $tp(\bar{b}/\bar{a}), tp(\bar{a}/\bar{d})$, and so $tp(\bar{a}/\bar{c})$ are semi-isolated. Thus by Fact 2.1, $tp(\bar{b}/\bar{c})$ is semi-isolated, while $tp(\bar{b}\bar{c}) = tp(\bar{b}\bar{e}) = tp(\bar{a}\bar{b})$, a contradiction.

# 3 Proof of Theorem 1.2

In this section, $T$ will be a countable, non $\aleph_0$-categorical theory.

**Fact 3.1** *(folklore) Suppose that $T$ has finitely many nonisomorphic models. Then there is a tuple $\bar{a}$ and a prime model $M$ over $\bar{a}$ such that $tp(\bar{a})$ is nonisolated and every complete n-type (for all n) over $\phi$ is realized in $M$. Moreover there is a tuple $\bar{b}$ in $M$ such that, $tp(\bar{b}) = tp(\bar{a})$ and $tp(\bar{a}/\bar{b})$ is nonisolated.*

*Proof.* Let $q_0, q_1.q_2, ...$ be an enumeration of all complete types of $T$ over $\phi$. Suppose that $\bar{e}_i \models q_i$ and $\bar{d}_i = \bar{e}_0\bar{e}_1...\bar{e}_i$. Now there is a prime model $N_i$ over $\bar{d}_i$ for each $i < \omega$. Thus for some $j < \omega$, $N_j(= M)$ is isomorphic to $N_i$ for infinitely many $i \geq j$. Therefore the prime model $M$ over $\bar{d}_j(= \bar{a})$ realizes every complete types over $\phi$. As $M$ is not prime over $\phi$, $tp(\bar{a})$ is not isolated.

Now since $T(\bar{a})$ is again non $\aleph_0$-categorical, for some tuple $\bar{s}$, $tp(\bar{s}/\bar{a})$ is nonisolated. Let $\bar{s}'\bar{b}(\in M)$ realize $tp(\bar{s}\bar{a})$. Then as $tp(\bar{s}'/\bar{b})$ is nonisolated, $M$ is not prime over $\bar{b}$. Since $M$ is prime over $\bar{a}$, $tp(\bar{a}/\bar{b})$ must not be isolated. □

We are ready to prove Theorem 1.2. We will use the same notation in the preceding Fact 3.1. Let $T$ be supersimple, and have finitely many models. We will lead a contradiction.

**Claim 3.2** Let $p = tp(\bar{a})$. There are two realizations $\bar{a}_0, \bar{a}_1$ of $p$ such that $\{\bar{a}_0, \bar{a}_1\}$ is independent (over $\phi$), and $tp(\bar{a}_0/\bar{a}_1)$ is nonisolated.

*Proof.* Let $\bar{c}$ be a realization of $p$ such that $tp(\bar{c}/\bar{a}\bar{b})$ does not fork over $\phi$. Now, by Fact 2.2, $tp(\bar{a}/\bar{b})$ is nonsemi-isolated. Hence, by Fact 2.1, either $tp(\bar{a}/\bar{c})$ or $tp(\bar{c}/\bar{b})$ must not be isolated. Thus $\bar{a}, \bar{c}$ or $\bar{c}, \bar{b}$ are desired two realizations of $p$. □

Now in the preceding claim, we may assume $\bar{a}_0, \bar{a}_1$ are in $M$. Moreover, as $\{\bar{a}_0, \bar{a}_1\}$ is independent, $tp(\bar{a}_1/\bar{a}_0)$ is also nonisolated, by Corollary 2.4. Now then $tp(\bar{a}/\bar{a}_0)$, $tp(\bar{a}/\bar{a}_1)$ are both nonisolated; for example if $tp(\bar{a}/\bar{a}_0)$ were isolated, then $M$ is prime over $\bar{a}_0$ and so $tp(\bar{a}_1/\bar{a}_0)$ were isolated, a contradiction. Therefore again by Corollary 2.4, $tp(\bar{a}/\bar{a}_0)$ and $tp(\bar{a}/\bar{a}_1)$ both fork over $\phi$.

Let us here summarize the relationships between three realizations $\{\bar{a}, \bar{a}_0, \bar{a}_1\}$ of $p$.

(1) $\{\bar{a}_0, \bar{a}_1\}$ is independent.

(2) For each $i = 0, 1$, $tp(\bar{a}_i/\bar{a})$ is isolated, whereas $tp(\bar{a}/\bar{a}_i)$ is nonisolated (so nonsemi-isolated). Thus $\{\bar{a}, \bar{a}_i\}$ is not independent.

Now then we are able to construct a tree $\{\bar{a}_\sigma | \sigma \in 2^{<\omega}\}$ such that $\bar{a}_\phi = \bar{a}$ and $tp(\bar{a}_\sigma \bar{a}_{\sigma 0} \bar{a}_{\sigma 1}) = tp(\bar{a} \bar{a}_0 \bar{a}_1)$ for each $\sigma \in 2^{<\omega}$ (**). Moreover, the basic properties of nonforking together with (1) enable us to assume that every antichain in the tree is independent, (e.g. $\{\bar{a}_{\bar{0}1} : |\bar{0}| = n \text{ for some } n < \omega\}$ is independent). Now by (2) with Fact 2.1, for each $\sigma \in 2^{<\omega}$ and each $i = 0, 1$, $tp(\bar{a}_{\sigma i}/\bar{a})$ is semi-isolated. But $tp(\bar{a}/\bar{a}_{\sigma i})$ is nonsemi-isolated, since if it were, then again by Fact 2.1, $tp(\bar{a}_\sigma/\bar{a}_{\sigma i})$ is semi-isolated, contradicting (2) and (**). Hence by Proposition 2.3, $tp(\bar{a}/\bar{a}_{\sigma i})$ forks over $\phi$.

Conclusively, we have countably many independent realizations of $p$, each of which is not independent with $\bar{a}$. Finally, by the symmetry and transitivity of nonforking, there is a sequence of complete types $\langle p_k | k \in \omega \rangle$ such that $p_0 = p$ and $p_{k+1}$ is a forking extension of $p_k$ for each $k \in \omega$. This violates supersimplicity of $T$. Therefore Theorem 1.2 is proved.